\documentclass[12pt]{article} 

\usepackage[leqno]{amsmath}
\usepackage{amssymb}
\usepackage{amsthm}
\usepackage{pstricks}
\usepackage{enumerate}
\newtheorem{theorem}{Theorem}
\newtheorem{lemma}{Lemma}

\newcommand{\mo}{\:\operatorname{mod}\:}
\newcommand{\ph}{\varphi}

\newcommand{\eps}{\varepsilon}

\newcommand{\zeile}{\vspace{\baselineskip}}

\newcommand{\N}{\mathbb{N}}

\newcommand{\R}{\mathbb{R}}

\allowdisplaybreaks

\begin{document}

\begin{center}
\begin{large}
\textsc{On the ternary Goldbach problem with primes in arithmetic
  progressions of a common module}
\end{large}
\end{center}
\begin{center}
Karin Halupczok
\end{center}

\begin{abstract}\noindent
For $A,\eps>0$ and any sufficiently large odd $n$
we show that for almost all $k\leq R:=n^{1/5-\eps}$
there exists a representation $n=p_{1}+p_{2}+p_{3}$ with
primes $p_{i}\equiv b_{i}$ mod $k$ for almost all admissible
triplets $b_{1},b_{2},b_{3}$ of reduced residues mod $k$.
\end{abstract}


\subsection{Introduction and results}

Let $n$ be a sufficiently large integer, consider an integer $k$
and let $b_{1},b_{2},b_{3}$ be integers that are relatively 
prime to $k\geq 1$, we assume that $0\leq b_{i}<k$, $i=1,2,3$.

We consider the ternary Goldbach problem of writing $n$ as
\begin{equation*}
  n=p_{1}+p_{2}+p_{3}
\end{equation*}
with primes $p_{1}$, $p_{2}$ and $p_{3}$ satisfying the 
three congruences
\begin{equation*}
   p_{i}\equiv b_{i}\mo k, \; i=1,2,3
\end{equation*}
for the common module $k$.
One is interested in the solvability of this question 
for all sufficiently large $n$ with the module $k$ 
being as large as up to some power of $n$.
This problem has been studied intensely by many authors.
For an overview, see for example \cite{c1}.

A necessary condition for solvability is
\begin{equation*}
  n\equiv b_{1}+b_{2}+b_{3} \mo k,
\end{equation*}
otherwise no such representation of $n$ is possible.

We call such a triplet $b_{1},b_{2},b_{3}$
of reduced residues mod $k$ \textit{admissible},
and a pair $b_{1},b_{2}$ of reduced residues \textit{admissible},
if $(n-b_{1}-b_{2},k)=1$. For a given $b_{1}$ we call $b_{2}$
\textit{admissible}, if $b_{1},b_{2}$ is an admissible pair.
Let us denote the number of these admissible pairs respectively 
triplets by $A(k)$.

We precise our consideration of this strengthened ternary 
Goldbach problem in the following way. Let
\begin{equation*}
  J_{3}(n):= J_{k,b_{1},b_{2},b_{3}}(n) :=
      \sum_{\substack{ m_{1},m_{2},m_{3}\leq n\\ m_{1}+m_{2}+m_{3}=n \\
                              m_{i}\equiv b_{i}\: (k),\\i=1,2,3}} 
    \Lambda (m_{1})\, \Lambda (m_{2})\, \Lambda (m_{3}),
\end{equation*}
where $\Lambda$ is von Mangoldt's function. $J_{3}(n)$ goes closely
with the number of representations of $n$ in the way mentioned.

In this paper we prove that the deviation of $J_{3}(n)$
from its expected main term is
uniformly small for large moduli, namely in the following sense.
\zeile
\begin{theorem}
\label{Th1}
   For every $A,\eps>0$, every sufficiently large $n$ and 
   for $D\leq n^{1/5-\eps}$ it holds that
   \begin{equation*}  
   \mathcal{E}:=
     \sum_{D<k\leq 2D} \frac{k}{\ph(k)} \sum_{(b_{1},k)=1}
      \frac{1}{\ph(k)}\sum_{b_{2}\text{ adm.}}
          \left|  J_{3}(n)   - \frac{n^{2}}{k^{2}}\mathcal{S}(n,k)
           \right|  \ll\; \frac{n^{2}}{(\log n)^{A}}.
   \end{equation*} 
\end{theorem}

Here $\mathcal{S}(n,k)$ denotes the singular series for this
special Goldbach problem and depends on $b_{1}$, $b_{2}$ and $k$ likewise
$J_{3}(n)$ does; residue $b_{3}$ is simply $b_{3}\equiv n-b_{1}-b_{2} \;(k)$.
Namely, see \cite{c2}, for odd $n$ we have
\begin{equation*}
   \mathcal{S}(n,k)= C(k)\prod_{p|k}
   \frac{p^{3}}{(p-1)^{3}+1} \prod_{\substack{p|n\\p\nmid k}}
   \frac{(p-1)((p-1)^{2}-1)}{(p-1)^{3}+1} \prod_{p>2}
   \biggl( 1+ \frac{1}{(p-1)^{3}} \biggr),
\end{equation*}
where $p>2$ throughout, $C(k)=2$ for odd $k$ and $C(k)=8$ for even $k$.

As a consequence of Theorem \ref{Th1}, we prove in section
\ref{seceins} the following result.

\begin{theorem}
\label{th2}
Let $A,\eps>0$ and let $n\in\N$ be odd and sufficiently large.
Then for all $k\leq R:=n^{1/5-\eps}$ with at most
$\ll R(\log n)^{-A}$ many exceptions of them there exists a representation
$n=p_{1}+p_{2}+p_{3}$ with primes $p_{i}\equiv b_{i} \:(k)$
for all but $\ll A(k)(\log n)^{-A}$ many admissible triplets 
$b_{1},b_{2},b_{3}$. 
\end{theorem}

So there are few exceptions for $k$, and also the number of 
exceptions of admissible triplets is small
compared with the number $A(k)$ of all admissible triplets.

Let us compare this Theorem \ref{th2} with the result of J.\ Liu and
T.\ Zhang in \cite{c2} who show the assertion for $R:=n^{1/8-\eps}$
and \textit{all} admissible triplets. In another paper \cite{c4}, Z.\
Cui improved this to $R:=n^{1/6-\eps}$. 
Further C.\ Bauer and Y.\ Wang showed in \cite{c1}
the assertion for $R:=n^{5/48-\eps}$, but with only $\ll (\log n)^{B}$ 
many exceptions.
 
Here we improved the bound for $R$ again, but at the cost
of possible but few exceptions of admissible triplets.

\subsection{Proof of Theorem \ref{th2}}
\label{seceins}

First of all we give a lower bound for $A(k)$:
\begin{lemma}
\label{l1}
  For odd $n$ we have $A(k)\gg\frac{\ph(k)^{2}}{(\log k)^{3}}$, more
  accurate, for every reduced residue $b_{1}$ mod $k$ there are $\gg
  \frac{\ph(k)}{(\log k)^{3}}$ many reduced residues $b_{2}$ mod $k$
  with $(n-b_{1}-b_{2},k)=1$.
\end{lemma}

\textbf{Proof.} Fix a reduced residue $b_{1}$ mod $k$.
Now count the $b_{2}$ with $(b_{2},k)=(n-b_{1}-b_{2},k)=1$.
So $b_{2}$ is to choose such that for all prime divisors $p>2$ of $k$
we have $b_{2}\not\equiv 0\;(p)$ and $b_{2}\not\equiv n-b_{1}\;(p)$,
what makes $\geq p-2$ many possibilities for $b_{2}$ mod $p$, and
$\geq p^{l-1}(p-2)$ many possiblilities for $b_{2}$ mod $p^{l}$.
If $p=2$ for even $k$ we have an odd $b_{1}$, so $n-b_{1}$ is even
and therefore one can take $b_{2}\equiv 1 (2)$, so there are
$2^{\nu_{2}(k)-1}$ many possibilities for $b_{2}$ mod $2^{\nu_{2}(k)}$,
if $2^{\nu_{2}(k)}||k$.

Therefore the number of $b_{2}$ is at least
\begin{equation*}
   2^{\max\{0,\nu_{2}(k)-1\}} \prod_{\substack{p^{l}||k \\ p\neq 2}}
   p^{l-1}(p-2)= \ph(k)\prod_{\substack{p|k \\ p\neq 2}}\frac{p-2}{p-1}
\end{equation*}
with
\begin{align*}
  \prod_{\substack{p|k\\p\neq 2}}\;&\frac{p-1}{p-2} =
   \prod_{\substack{p|k\\p\neq 2}}
  \biggl(1+\frac{1}{p-2}\biggr) \leq \prod_{p|k}
  \biggl(1+\frac{2}{p-1}\biggr) \\
  &\leq \sum_{q=1}^{k} \frac{\mu(q)^{2}\, 2^{\omega(q)}}{\ph(q)}
  \ll \sum_{q=1}^{k}\frac{\tau(q)}{q}\log k\ll (\log k)^{3}.
\end{align*}
This shows the Lemma.  \hfill$\square$

Now we show Theorem \ref{th2} as a corollary of Theorem \ref{Th1}.

Fix $A,\varepsilon>0$ and let $n$ be odd and sufficiently large.
Consider 
\begin{equation*}
  R_{3}(n):=\sum_{\substack{p_{1},p_{2},p_{3}\\p_{1}+p_{2}+p_{3}=n 
  \\p_{i}\equiv b_{i} (k),\\ i=1,2,3 }}  \log p_{1}\log p_{2} \log p_{3}
\quad\text{ and }\quad
  r_{3}(n):=\sum_{\substack{p_{1},p_{2},p_{3}\\p_{1}+p_{2}+p_{3}=n 
  \\p_{i}\equiv b_{i} (k),\\ i=1,2,3 }} 1.
\end{equation*}

Let $D<k\leq 2D$ with $D\leq R:=n^{1/5-\eps}$. For any
admissible triplet $b_{1},b_{2},b_{3}$ mod $k$ we have
\begin{equation*}
  \left| R_{3}(n) - J_{3}(n) \right| \leq (\log n)^{3} W_{3},
\end{equation*}
where $W_{3}$ denotes the number of solutions of $p^{l}+q^{j}+r^{m}=n$
with $p,q,r$ prime, and where $l,j$ or $m$ are at least 
$2$ such that $p^{l}\equiv b_{1}\:(k)$,
$q^{j}\equiv b_{2}\:(k)$ and $r^{m}\equiv b_{3}\:(k)$.

Now we prove that 
\begin{equation*}
\sum_{D<k\leq 2D} k
\max_{\substack{b_{1},b_{2},b_{3}\\\text{admissible}}} W_{3}
   \ll \frac{n^{2}}{(\log n)^{A+3}}.
\end{equation*}

For this, we split the number $W_{3}$ according to if at least two of the
exponents $l,j,m$ are $\geq 2$ or only one, and for this
we write $W_{3}=W_{1}+W_{2}$. There are at most
$\sqrt{n}$ prime powers $\leq n$ with exponent $\geq 2$,
so in the first case we have $W_{1}\ll n$, and the left hand 
side with $W_{1}$ is $\ll D^{2}W_{1}\ll
D^{2}n\ll\frac{n^{2}}{(\log n)^{A+3}}$.

In the second case, if only one exponent is $\geq 2$,
we have $W_{2}\ll \sqrt{n}\cdot\frac{n}{k}=\frac{n^{3/2}}{k}$, and so
the left hand side is  $\ll Dn^{3/2}\ll \frac{n^{2}}{(\log n)^{A+3}}$.

So for $D\leq n^{1/5-\eps}$ it follows from Theorem \ref{Th1}:

\begin{align*}
 &\sum_{D<k\leq 2D}
 \frac{k}{\ph(k)} \sum_{(b_{1},k)=1} 
   \frac{1}{\ph(k)} \sum_{b_{2}\text{ adm.}}
 \biggl| R_{3}(n) - \frac{n^{2}}{k^{2}}\mathcal{S}(n,k) \biggr| \\ 
 \leq &\sum_{D<k\leq 2D} \frac{k}{\ph(k)} \sum_{(b_{1},k)=1} 
   \frac{1}{\ph(k)} \sum_{b_{2}\text{ adm.}}
   \biggl| R_{3}(n)- J_{3}(n) \biggr| \\&+  
   \sum_{D<k\leq 2D} 
 \frac{k}{\ph(k)} \sum_{(b_{1},k)=1} 
   \frac{1}{\ph(k)} \sum_{b_{2}\text{ adm.}}
   \biggl|J_{3}(n) - \frac{n^{2}}{k^{2}}\mathcal{S}(n,k) \biggr| \\ 
  \ll &(\log n)^{3} \sum_{D<k\leq 2D}  k
  \max_{\substack{b_{1},b_{2},b_{3}\\\text{admissible}}} W_{3}  +
  \frac{n^{2}}{(\log n)^{A}}\ll  \frac{n^{2}}{(\log n)^{A}}.
\end{align*}
So the formula of Theorem \ref{Th1} holds also for $R_{3}(n)$ instead
of $J_{3}(n)$. 

Now for $D<k\leq 2D$
we have
 $A(k):=\#\{b_{1},b_{2}\text{ admissible mod }k\}$, and let \\
$T(k):=\#\{b_{1},b_{2}\text{ admissible mod }k; \; R_{3}(n)=0\}$
and consider the set
\begin{equation*}
   \mathcal{K}_{D}:=\{k; \; D<k\leq 2D, T(k)\geq A(k)(\log n)^{-A} \}
\end{equation*}
and let $K_{D}$ be its number.

Since $\mathcal{S}(n,k)\gg 1$ if it is positive, 
what is the case for admissible triplets and odd $n$
(see its formula above as an Euler product), we have
\begin{align*}
  K_{D}\cdot \frac{n^{2}}{D} &\ll \sum_{\substack{D<k\leq 2D \\ k\in
      \mathcal{K}_{D}}} \frac{k}{T(k)}
  \sum_{\substack{b_{1},b_{2}\text{ adm.} \\ R_{3}(n)=0} } \;\biggl|
   \frac{n^{2}}{k^{2}}\mathcal{S}(n,k)\biggr| \\
    &\ll \sum_{D<k\leq 2D} \frac{k}{A(k)} \sum_{b_{1},b_{2}\text{ adm.}} (\log
    n)^{A} \;\biggl|R_{3}(n)-\frac{n^{2}}{k^{2}}\mathcal{S}(n,k)\biggr| \\
    &\ll (\log n)^{A+3} \sum_{D<k\leq 2D} \frac{k}{\ph(k)^{2}}
    \sum_{\substack{b_{1},b_{2} \\ \text{adm.}}}
     \;\biggl|R_{3}(n)-\frac{n^{2}}{k^{2}}\mathcal{S}(n,k)\biggr|  
   \ll \frac{n^{2}}{(\log n)^{A}}, 
\end{align*}
using Lemma \ref{l1} and the above.
Therefore it follows that $K_{D}\ll D(\log n)^{-A}$, so
for all $k\not\in \mathcal{K}_{D}$ we have 
$R_{3}(n)>0$ for all but $\ll A(k)(\log n)^{-A}$ many
admissible triplets $b_{1},b_{2},b_{3}$,
and then $r_{3}(n)\gg R_{3}(n)(\log n)^{-3}$ is positive, too.
This shows Theorem \ref{th2}, since the overall number of exceptions is
\begin{equation*}
  \ll\sum_{i\ll \log R} \#\mathcal{K}_{2^{i}} \ll (\log
  n)\cdot\frac{R}{(\log n)^{A+1}}=\frac{R}{(\log  n)^{A}}.
\end{equation*} \hfill$\square$

\subsection{Proof of Theorem \ref{Th1}}
\label{secneu}

We are going to show Theorem \ref{Th1} in two steps
according to the circle method.

Let $A,\eps,\theta>0$, $B\geq 2A+1$ and $D\leq n^{1/4}(\log n)^{-\theta}$.

We define major arcs $\mathfrak{M}\subseteq\mathbb{R}$ by
\begin{equation*}
  \mathfrak{M}:= \bigcup_{q\leq D(\log n)^{B}}
     \bigcup_{\substack{0<a<q \\ (a,q)=1}}
  \left]\frac{a}{q}-\frac{D(\log n)^{B}}{qn}, 
     \frac{a}{q}+\frac{D(\log n)^{B}}{qn} \right[
\end{equation*}
and minor arcs by
\begin{equation*}
\mathfrak{m}:=\left]-\frac{D(\log n)^{B}}{n},1-\frac{D(\log n)^{B}}{n}\right[
\;\backslash\;\mathfrak{M}.
\end{equation*}

For $\alpha \in \mathbb{R}$ and some residue $b$ mod $k$ denote
\begin{equation*}
S_{b}(\alpha) := S_{b,k} (\alpha)
           := \sum_{ \substack{m\leq n\\m\equiv b\,(k)}}
                    \Lambda (m)  \; e(\alpha m).
\end{equation*}

From the orthogonal relations for $e(\alpha m)$ it follows that
\begin{equation*}
  J_{3}(n) = \int_{0}^{1} S_{b_{1}}(\alpha) 
                S_{b_{2}}(\alpha) S_{b_{3}}(\alpha)\,
                 e(-n\alpha)\, d\alpha.
\end{equation*}

By
\begin{equation*}
  J^{\mathfrak{M}}_{3} (n):= \int_{\mathfrak{M}} S_{b_{1}}(\alpha) 
                S_{b_{2}}(\alpha) S_{b_{3}}(\alpha)
                 \,e(-n\alpha)\, d\alpha
\end{equation*}
we denote the value of the integral for $J_{3}(n)$ on the major arcs
$\mathfrak{M}$ and by
\begin{equation*}
  J_{3}^{\mathfrak{m}}(n):= J_{3} (n) - J_{3}^{\mathfrak{M}} (n)
\end{equation*}
the value on the minor arcs $\mathfrak{m}$. 

Concerning the major arcs we have
\begin{theorem} 
\label{Th2}
For $D \leq n^{1/5-\eps}$ it holds that
  \begin{equation*}
    \mathcal{E}^{\mathfrak{M}}:=
      \sum_{D<k\leq 2D} k
      \max_{\substack{b_{1},b_{2},b_{3}\\\text{admissible}}}
      \left|J_{3}^{\mathfrak{M}} (n) -
      \frac{n^{2}}{k^{2}}\mathcal{S}(n,k)\right|
      \ll \frac{n^{2}}{(\log n)^{A}}.
  \end{equation*}
\end{theorem}

We can give the proof of Theorem \ref{Th2} very shortly, as
it is simply done by adapting the result of J.\ Liu and T.\ Zhang in
\cite{c2} for the here given major arcs. In fact, by pursuing their proof
we see that for $P:= D(\log n)^{B}$ and $Q:=\frac{n}{D(\log n)^{B}} $
and any $U\leq P$, we have to choose $D$ such that the conditions
\begin{xalignat*}{2}
  U&\leq n^{1/3}(\log n)^{-E}, & (UQ)^{-1}&\leq U^{-3}(\log n)^{-E}\\
  DU&\leq D^{1/3-\delta}n^{1/3}(\log n)^{-E}, & (UQ)^{-1}&\leq
  D^{1-\delta}(DU)^{-3} (\log n)^{-E}
\end{xalignat*}
are satisfied for any $E>0$ and small $\delta >0$.
The optimal choice of $D$ is therefore given by $D\leq n^{1/5-\eps}$
what proves Theorem \ref{Th2}. The improvement in this paper comes
from the different intervals given as major and minor arcs 
such that dealing on the minor arcs
with mean values over $b_{1},b_{2}$ is still possible.

Namely, as estimate on the minor arcs we show in the next section \ref{secdrei}:
\begin{theorem} 
\label{Th4}
For $D \leq n^{1/4}(\log n)^{-\theta}$ we have
  \begin{equation*}
    \mathcal{E}^{\mathfrak{m}}:=
      \sum_{D<k\leq 2D} \frac{k}{\ph(k)}
      \sum_{(b_{1},k)=1} \frac{1}{\ph(k)}
       \sum_{\substack{(b_{2},k)=1\\\text{adm.}}} 
       \left|J_{3}^{\mathfrak{m}} (n)\right|
      \ll \frac{n^{2}}{(\log n)^{A}}.
  \end{equation*}
\end{theorem}

Theorem \ref{Th1} is then a corollary of Theorems \ref{Th2} and
 \ref{Th4} since $\mathcal{E}\leq
\mathcal{E}^{\mathfrak{M}}+\mathcal{E}^{\mathfrak{m}}$.

This Theorem \ref{Th4} is the interesting part of Theorem \ref{Th1},
where we can gain a higher power of $n$ for the bound of $D$ by considering the
mean value over $b_{1},b_{2}$ instead of the maximum. But due to this
we have to allow exceptions of admissible triplets in Theorem
\ref{th2}, as we have seen in its proof.

In both Theorems \ref{Th2} and \ref{Th4} the resulting bound for $D$
is the optimum with the given method, these bounds cannot be
balanced to get a larger range than $n^{1/5}$.
Also the cited method for the major arcs cannot be improved to
gain from mean values over $b_{1},b_{2}$ since the used
character sum estimates are independent of $b_{1},b_{2}$.
But it may be possible that another method would succeed on
$\mathfrak{M}$.

\subsection{Proof of Theorem \ref{Th4}, the estimate on the minor arcs}
\label{secdrei}


Let $D \leq n^{1/4}(\log n)^{-\theta}$
and consider $\mathcal{E}^{\mathfrak{m}}$, it is (where
$b_{1},b_{2}$ run through all reduced residues mod $k$ if indicated by
a star)
\begin{align*}
 \ll &\sum_{D<k\leq 2D} \frac{k}{\ph(k)^{2}}  \sideset{}{^{*}}\sum_{b_{1},b_{2}}
      \left|J_{3}^{\mathfrak{m}} (n)\right| \\
  \leq &\sum_{D<k\leq 2D} \frac{k}{\ph(k)^{2}} \sideset{}{^{*}} 
     \sum_{b_{1},b_{2}} \int_{\mathfrak{m}}
     \left|S_{b_{1}}(\alpha)S_{b_{2}}(\alpha)S_{n-b_{1}-b_{2}}(\alpha)
     \right| d\alpha \\
  = &\sum_{D<k\leq 2D} \frac{k}{\ph(k)}\sideset{}{^{*}}\sum_{b_{1}} \int_{\mathfrak{m}}
  |S_{b_{1}}(\alpha)| \cdot \frac{1}{\ph(k)}\sideset{}{^{*}}\sum_{b_{2}}
  |S_{b_{2}}(\alpha)S_{n-b_{1}-b_{2}}(\alpha)|d\alpha \\
  \leq &\sum_{D<k\leq 2D} \frac{k}{\ph(k)} \int_{\mathfrak{m}}
  \sideset{}{^{*}}\sum_{b_{1}} |S_{b_{1}}(\alpha)| \\
       &\qquad\cdot \frac{1}{\ph(k)} \biggl( \sum_{b_{2}\text{ mod } k}
     |S_{b_{2}}(\alpha)|^{2}  \biggr)^{1/2} \biggl(\sum_{b_{2} \text{ mod } k}
     |S_{n-b_{1}-b_{2}}(\alpha)|^{2}  \biggr)^{1/2} d\alpha \\
  \leq & \sum_{D<k\leq 2D} \frac{k}{\ph(k)}  \max_{\alpha\in\mathfrak{m}}
     \sideset{}{^{*}}\sum_{b_{1}} |S_{b_{1}}(\alpha)|  \frac{1}{\ph(k)}
     \sum_{b_{2}\text{ mod } k} \int_{0}^{1}|S_{b_{2}}(\beta)|^{2} d\beta   \\ 
  \ll & n(\log n)^{3}\sum_{D<k\leq 2D} \frac{1}{\ph(k)} \max_{\alpha\in\mathfrak{m}}
    \sideset{}{^{*}} \sum_{b_{1}} |S_{b_{1}}(\alpha)| \\ \leq & n
    (\log n)^{3} \sum_{D<k\leq 2D} 
       \max_{\alpha\in\mathfrak{m}}
       \biggl(\frac{1}{\ph(k)}\sideset{}{^{*}}\sum_{b_{1}} 
     |S_{b_{1}}(\alpha)|^{2} \biggr)^{1/2} \\
  \ll &   n (\log n)^{3} \sum_{D<k\leq 2D}
    \biggl(\frac{n^{2}}{D^{2}(\log n)^{2A+6}}\biggr)^{1/2}.
\end{align*}
In the last step we use Lemma \ref{L2} that will be shown next,
valid for $D\leq n^{1/4}(\log n)^{-\theta}$ and suitable chosen
$\theta,B>0$ depending just on $A>0$.

Now  the above is $\ll n(\log n)^{3}D\frac{n}{D(\log n)^{A+3}}\ll
\frac{n^{2}}{(\log n)^{A}}$ as was to be shown for the minor arcs.
\hfill $\square$

So what is left is to show 
\begin{lemma}
   \label{L2} 
   For all $A>0$ and $B\geq 2A+1,\theta\geq B/2$
   let $D\leq n^{1/4}(\log n)^{-\theta}$ and $\alpha\in\R$ with
   $||\alpha-\frac{u}{v}||<\frac{1}{v^{2}}$ for some integers $u,v$
   with $(u,v)=1$ and $D(\log n)^{B}\leq v \leq \frac{n}{D(\log
     n)^{B}}$. Then for $D<d\leq 2D$ we have
   \[
      \frac{1}{\ph(d)}\sum_{c,(c,d)=1} |S_{c,d}(\alpha)|^{2} \ll
      \frac{n^{2}}{D^{2}(\log n)^{A}}.
   \]
\end{lemma}
 
We remark that for $\alpha\in\mathfrak{m}$ there exist $u,v$ with
$(u,v)=1$, $v\leq \frac{n}{D(\log n)^{B}}$ and 
$||\alpha-\frac{u}{v}||<\frac{D(\log n)^{B}}{vn}\leq\frac{1}{v^{2}}$
by Dirichlet's approximation theorem, so $v\geq D(\log n)^{B}$ since
$\alpha\in\mathfrak{m}$, and therefore the conditions of Lemma \ref{L2} are
fulfilled.

For the proof we need the following well known auxiliary Lemma.
\begin{lemma}
   \label{aux}
   Let $||\alpha-\frac{u}{v}||\leq\frac{1}{v^{2}}$, $(u,v)=1$. Then
   \[
      \sum_{m\leq X} \min(Y,||\alpha m||^{-1}) \ll \frac{XY}{v} +
      (X+v)(\log v).      
   \]  
\end{lemma}

\textbf{Proof of Lemma \ref{L2}.} 
Fix $n$ large and $D\leq n^{1/4}(\log n)^{-\theta}$,
and let $\alpha$, $u$ and $v$ be as given in Lemma \ref{L2}.

We apply Vaughan's identity on the exponential sum $S_{c,d}(\alpha)$,
see for example A.\ Balog in \cite{c3}, where a similar Lemma is shown
(Lemma 2 there). 
From that it follows that it suffices to show for any
complex coefficients $|a_{m}|,|b_{k}|\leq 1$ and any $M\in\N$ with
\begin{align*}
  &I: \qquad M\leq V^{2}, \text{ if } b_{k}=1 \text{ for all }k,\\
  &II: \qquad V\leq M\leq \frac{n}{V} \text{ else, where }V:=D(\log n)^{B},
\end{align*}
we have
\[
   \sum_{(c,d)=1} \Biggl|\sum_{m\sim M} \sum_{\substack{k\leq
       n/m \\ km\equiv c (d)}} a_{m} b_{k} e(\alpha mk)\Biggr|^{2}\ll
   \frac{n^{2}}{D(\log n)^{A}}.
\]

Here $m\sim M$ means $M<m\leq M'$ for some $M'\leq 2M$.

We consider first \textbf{case II}:
Then the left hand side becomes (where $m^{*}$ denotes the inverse of $m$
mod $d$):
\begin{align*}
  II&:=\sum_{(c,d)=1} \biggl|\sum_{\substack{m\sim M\\(m,d)=1}}
    a_{m} \sum_{\substack{k\leq n/m \\ km\equiv c (d)}} b_{k}
    \,e(\alpha mk)\biggr|^{2} \\
   &\leq \sum_{(c,d)=1} M \sum_{\substack{m\sim M\\(m,d)=1}}   
    \biggl|\sum_{\substack{k\leq n/m \\ k\equiv cm^{*} (d)}} b_{k} 
    \,e(\alpha mk)\biggr|^{2} \\
   &= M\sum_{\substack{m\sim M\\(m,d)=1}}
   \sum_{(c,d)=1} \biggl|\sum_{\substack{k\leq n/m \\ k\equiv c (d)}}
    b_{k}\, e(\alpha mk)\biggr|^{2} \\
   &= M \sum_{m\sim M} \sum_{(c,d)=1}
   \sum_{\substack{k\leq n/m \\ k\equiv c (d)}} b_{k}
   \sum_{\substack{k'\leq n/m \\ k'\equiv k (d)} } \overline{b_{k'}}
   \,e(\alpha m(k-k'))  \\
   &= M \sum_{m\sim M} \sum_{\substack{k\leq n/m\\(k,d)=1}} b_{k}
     \sum_{\substack{k'\leq n/m \\ k'\equiv k (d)}} \overline{b_{k'}}
     \,e(\alpha m(k-k'))  \\
   &=M \sum_{m\sim M} \sum_{\substack{k\leq n/m\\(k,d)=1}}
        b_{k} \sum_{\substack{l\geq (k-n/m)/d \\l\leq (n/m-1)/d}}
   \overline{b_{k-ld}}\, e(\alpha mld) \\
   &\leq M \sum_{k\leq n/M} \sum_{|l|\leq n/Md} 
      \biggl|\sum_{\substack{m\sim M \\ m\leq n/k \\ m\leq 
          n/\max\{k-ld,ld+1\}}} e(\alpha mld)\biggr|.
\end{align*}
Now the exponential sum in absolute value
is $\ll \min(M, ||\alpha ld||^{-1})$, so the
estimation goes on with
\begin{align*}
  &\ll M\frac{n}{M} \sum_{|l|\leq n/Md} \min(M,||\alpha ld||^{-1}) \\
  &\ll n \sum_{\substack{L\leq n/M\\d|L}} \min(M,||\alpha L||^{-1}) +nM\\
  &\leq n \biggl(\sum_{\substack{L\leq n/M\\d|L}}1^{2} \biggr)^{1/2}
   \biggl(\sum_{L\leq n/M}M\min(M,||\alpha L||^{-1}) \biggr)^{1/2} +nM\\
  &\ll n \biggl(\frac{n}{Md}\biggr)^{1/2} M^{1/2}
    \biggl( \frac{n}{v} + \biggl(\frac{n}{M}+v\biggr)(\log n) \biggr)^{1/2} +nM
\end{align*}
because of the auxiliary Lemma \ref{aux}. So expression $II$ is
$\ll\frac{n^{2}}{D(\log n)^{A}}$
since we have $D(\log n)^{B}=V\leq M\leq n/V$ in case II, and since
$D(\log n)^{B}\ll v\ll \frac{n}{D(\log n)^{B}}$ for $B\geq 2A+1$.

Now consider \textbf{case I}:
Then the left hand side becomes (where $m^{*}$ denotes the inverse of
$m$ mod $d$):
\begin{align*}
 I:= &\sum_{(c,d)=1} \biggl| \sum_{\substack{m\sim M\\(m,d)=1}}
    a_{m} \sum_{\substack{k\leq n/m \\ km\equiv c (d)}} e(\alpha
    mk)\biggr|^{2} \\ 
   \leq & \sum_{(c,d)=1}  M \sum_{\substack{m\sim M\\(m,d)=1}}
      \biggl|\sum_{\substack{k\leq n/m \\ k\equiv cm^{*} (d)}}
      e(\alpha mk)\biggr|^{2} \\
    \leq &M \sum_{m\sim M} \sum_{(c,d)=1} \biggl|\sum_{\substack{k\leq n/m
        \\ k\equiv c (d)}} e(\alpha mk) \biggr|^{2} \\ 
    = & M\sum_{m\sim M} \sum_{(c,d)=1} \sum_{\substack{k\leq n/m
       \\ k\equiv c (d)}} e(\alpha mk) \sum_{\substack{k'\leq n/m
       \\ k\equiv k' (d)}} e(-\alpha mk')  \\
    = & M \sum_{m\sim M} \sum_{\substack{k\leq n/m \\ (k,d)=1}}
     \sum_{\substack{k'\leq n/m \\ k\equiv k' (d)}} e(\alpha m (k-k')) \\
    \leq & M \sum_{m\sim M} \sum_{k\leq n/m} \biggl|
    \sum_{\substack{l\geq (k-n/m)/d \\ l\leq (n/m-1)/d}} e(\alpha mdl) \biggr| \\
    \ll & M \sum_{m\sim M} \sum_{k\leq n/M} \biggl(\min\biggl(\frac{n}{Md},||\alpha
    md||^{-1}\biggr) +1\biggr) \\
    \ll & n \sum_{m\sim M} \min\biggl(\frac{n}{Md},||\alpha
    md||^{-1}\biggr) + Mn \\
     \ll & n \sum_{\substack{L\sim Md \\ d|L}}
     \min\biggl(\frac{n}{Md},||\alpha L||^{-1}\biggr) + Mn \\
     \ll & n \biggl(\sum_{\substack{L\sim Md \\ d|L}}1\biggr)^{1/2}
          \biggl(\sum_{L\sim Md} \frac{n}{Md}
          \min\biggl(\frac{n}{Md},||\alpha L||^{-1}\biggr) \biggr)^{1/2} + Mn \\
     \ll & n \biggl(\frac{Md}{d}\biggr)^{1/2} 
       \biggl(\frac{n}{Md}\biggr)^{1/2} \biggl(\frac{n}{v} + (Md+v)(\log
       n)\biggr)^{1/2} + Mn 
\end{align*}
using again the auxiliary Lemma \ref{aux}. Now we get $I\ll
\frac{n^{2}}{D(\log n)^{A}} $ since
$D(\log n)^{B}\ll v \ll \frac{n}{D(\log n)^{B}}$ with $B\geq 2A+1$
and since $Md\ll V^{2}d\ll D^{3}(\log n)^{B} \ll \frac{n}{D(\log
  n)^{B}}$ for $D\leq n^{1/4}(\log n)^{-\theta}$ and $\theta\geq B/2$.
So Lemma \ref{L2} is shown.
\hfill $\square$

\vspace{\baselineskip}
\textbf{Remark added by author.} As was kindly pointed out
to me by Z. Cui, it is possible to improve
the statement on the major arcs such that
Theorems 1, 2 and 3 hold for the improved exponent
$1/4$ instead of $1/5$. This major arc improvement
has its idea in the publication of Z. Cui in \cite{c4}.



\end{document}